# 0-rotatability of classes of rooted symmetric trees. Are rooted symmetric trees 0-rotatable?

Rafael I. ROFA

**ABSRACT.** A graceful labelling of a tree T = (V,E), where V is the set of vertices of the tree and E is its edge set, is an bijective function f : V → {0,1,2, …, |E|} such that if edge uv is assigned the label g(uv) = |f(u) - f(v)| then the function g : E → {1, …, |E|} is also a bijection (that is all edge labels are distinct). A tree is said to be *0-roratable if for any of its vertices there is a graceful labelling that assigns the label 0 to that vertex*. A rooted symmetric tree *is a tree in which all vertices at the same level from root vertex have the same degree*. It was known since 1979 that rooted symmetric trees are graceful and an algebraic definition of graceful labelling of this class of trees was found by the author. In this paper we prove that rooted symmetric trees with at most 3 levels (including root vertex) are 0-rotatable. We also prove that symmetric spider trees with leg length at most 3 and symmetric banana trees, both of which are classes of rooted symmetric trees with 4 levels, are 0-rotatable. Based on these results, we conjecture that all spiders are 0-rotatable and raise the more general question whether all symmetric rooted trees are 0-rotatable.

## 1. Introduction

A graceful labelling of a tree T = (V, E) with vertex set V and edge set E, *is an injective function f : V → {0, 1, …, |E|} such that if edge uv is assigned the label g(uv)=|f(u)-f(v)| then the function g : E → {1, …, |E|}  is also injective* (that is all edge labels are distinct). Rosa (1967) first introduced several graph labelling methods as a means of attacking the Ringel Conjecture that the complete graph $K_{2n+1}$ can be decomposed into 2n+1 isomorphic trees each on n edges (Ringel, 1963) which is demonstrably true if one can prove that all trees are graceful, hence the birth of one of the most famous unsolved conjectures in graph theory, namely the graceful tree conjecture. Golomb (1972), coined the phrase 'graceful labellings' to refer to *β*-valuations which is one of the several graph labellings introduced initially by Rosa in 1967. Numerous classes of trees have been shown to be graceful and many of these are repeatedly listed in almost every paper about the subject (See for example, (J. C.  Bermond, 1979; Rofa, 2021; Superdock, 2014) A comprehensive list of graceful graphs and other labelling schemes of graphs is found in the electronic survey by J.A Gallian (Gallian, 2022) and in the survey by Robeva (2011). A tree T is said to be 0-rotatable *if there exists a graceful labelling that assigns the label 0 to any given vertex of the tree*. 0-rotatability of a tree is, by definition, a stronger condition than the existence of a graceful labelling for it, and hence a 0-rotatable tree is graceful but not every graceful tree is necessarily 0-rotatable and there are graceful trees that are not 0-rotatable (Bussel, 2004). The importance of the allocation of the label 0 in a graceful labelling of a tree stems from the fact that it plays a role in constructing larger graceful trees. For example, a gracefully labelled tree with n vertices can be gracefully 'grown' by adding k new leaves (vertices) with labels n, n+1, …, n+k to the vertex labelled 0. Moreover, it is also possible to combine a tree with an $\alpha$-labelling and any tree with a graceful labelling by identifying the vertices labelled 0 in a way such that the resultant tree is graceful. (An $\alpha$-labelling of a graph G is a graceful labelling f of G with the additional property that there exists an integer $k \in \{0, …, |EG|\}$ such that for each edge $uv \in EG, f(u) \leq k < f(v)$).

Rosa (1967) showed that path graphs are graceful and in 1977 he proved that this class of trees is 0-rotatable (Rosa, 1977). There are relatively much fewer results about the 0-rotatability of classes of trees than there are about classes of graceful trees and the 0-rotatability question is still open even for caterpillars which, together with numerous other classes of trees, have long been shown to be graceful. Chung and Hwang (1981) showed that caterpillars whose non-leaf vertices have the same degree are 0-rotatable

Bussel (2004) showed that trees with diameter at most 3 are 0-rotatable. He also showed that there are caterpillars with diameter 4 that are not 0-rotatable and that all non-0-rotatable trees with at most 14 vertices either are caterpillars with diameter four, or are trees formed by identifying the central vertex of a non-0-rotatable tree of diameter four with the end of a path $P_n$, $n \geq 1$ (which is not a caterpillar unless n=1 in which case it is a caterpillar with diameter 4). Consequently, Bussel conjectured that all non-0-rotatable trees belong to these two families of trees. It is noted that if Bussel's conjecture is true, then every caterpillar with diameter at least five is 0-rotatable. Luiz, Campos, and Richter (2016) proved that all caterpillars with diameter five or six are 0-rotatable. They also proved that if T is a caterpillar with diam(T) $\geq$ 7 and, for every non-leaf vertex $v \in V(T)$, the number of leaves adjacent to v is at least $2 + 2((\text{diam}(T) − 1) \mod 2)$, then T is 0-rotatable (Luiz, Campos, & Richter, 2017).

In this paper we investigate the 0-rotatability of rooted symmetric trees and some subclasses of rooted symmetric trees. A rooted symmetric tree is a tree with a root vertex and equal degree vertices at each given level of the tree, (see figure 1). A more precise definition is given in the next section. (Notice that rooted symmetric trees are mostly not caterpillars except for specific, trivial values of the parameters and hence are mostly uncovered by work done on the problem of 0-rotatability of caterpillars).

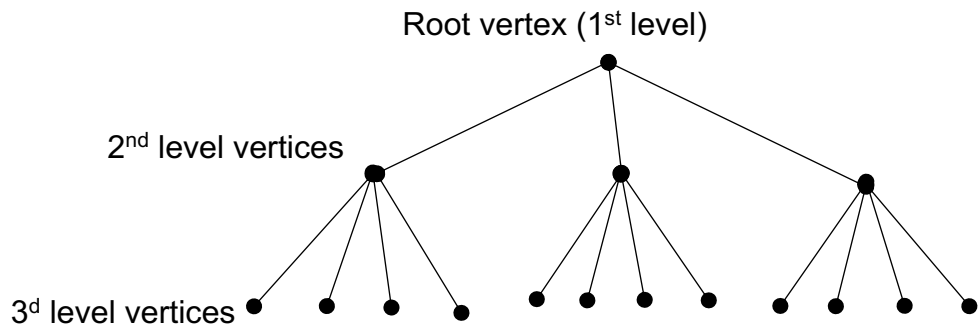

*Figure 1. A rooted symmetric tree*

J. C. Bermond and Sotteau (1976) proved that rooted symmetric trees are graceful. This was followed by many other alternative proofs (Robeva, 2011; Rofa, 2021; Sandy, Rizal, Manurung, & Sugeng, 2018).

In this paper, we prove the following results:

(1) Rooted symmetric trees with 3 levels including root vertex as the first level (diameter 4) are 0-rotatable. We also prove the following about subclasses of rooted symmetric trees of diameter 6 as our 2$^{nd}$ and 3$^d$ objectives:
(2) Symmetric spiders with leg length at most 3 are 0-rotatable. *A spider tree is a tree with at most one vertex of degree greater than* 2. The vertex of degree greater than 2, if it exists, is called the branch point of the tree. A leg of a spider tree is any path from the branch points to a leaf of the tree. We define a symmetric spider *to be a spider with equal length legs.* Note that symmetric spiders form a subclass of rooted symmetric trees (see figure 2). We propose the conjecture that all symmetric spiders are 0-rotatable.

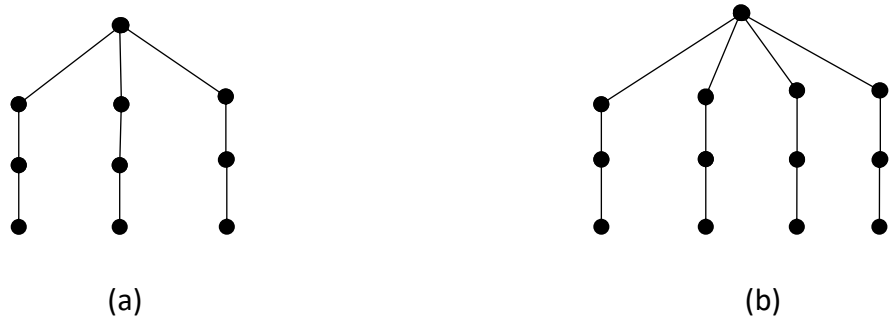

*Figure 2. Symmetric spiders with leg length 3 and 3 branches (a) and 4 branches (b)*

(3) Symmetric banana trees. A banana tree *consists of a vertex v joined to one leaf of any number of stars.* We define a symmetric banana tree *to be a banana tree that consists of a vertex v joined to one leaf of any number of the similar copies of a star.* (see figure 3 for an example of a symmetric banana tree).

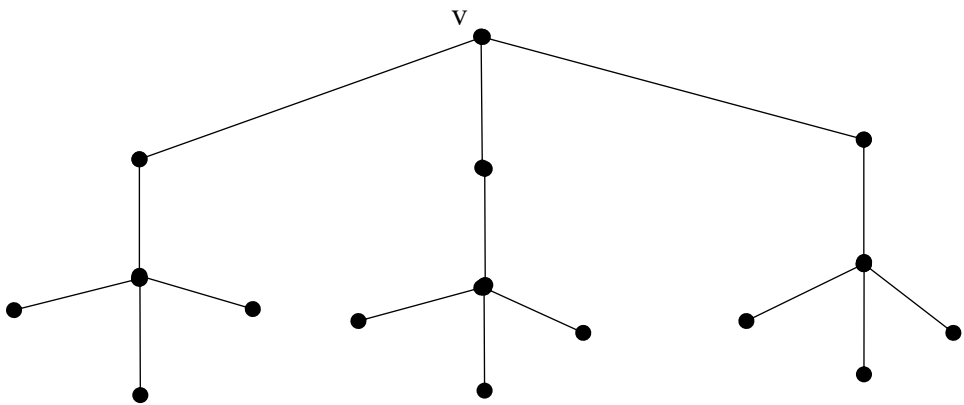

*Figure 3. A symmetric banana tree*

Finally, we raise the question: Are symmetric spiders with leg length 4 or more 0-rotatable? And generally: are rooted symmetric trees with 4 or more levels 0-rotatable?

## 2. Terminology and Preliminary results

As mentioned earlier there are many approaches that have been used by various authors to prove that rooted symmetric trees are graceful (J. C. Bermond & Sotteau, 1976; Robeva, 2011). However, before 2021 such proofs were either algorithmic or inductive in nature. In 2021, the author of this paper proved that rooted symmetric trees are graceful by constructing an algebraic graceful function that *maps* each vertex of a rooted symmetric tree to a unique number (Rofa, 2021). For the purpose of achieving the 3 goals of this paper, it is necessary to present here the terminology and the result (without proof) of the paper by the author Rofa (2021) as the same terminology and result will be used for achieving the goals of this paper. This done next.

let T be a rooted symmetric tree with root vertex $v_0$ (which is considered the only vertex at level 1). The daughter degree sequence of T is a sequence $S = (k_1, k_2, \ldots, k_{q-1})$, where $k_1$ is the degree of root vertex, $k_2$ is the degree of a vertex at level 2 (not counting the edge that connects the vertex to the preceding level 1), …., $k_{q-1}$ is the degree of vertex at level q - 1 (again not counting the edge that connects the vertex to preceding level) and $k_q = 1$ (no daughter vertices) is the degree of any vertex at the last level, which is level q. The edges in each branch at any level say level $i$ are indexed by labels $0, 1, \ldots, k_i - 1$ (not to be confused with the edges induced by the graceful labelling) and this is applied to all branches at this level in the same direction (from left to right). Thus, any vertex v of T at any level, say level r, is identified by a unique sequence of edge indices, say $(x_1, x_2, \ldots, x_{r-1})$. This setup is clarified for the rooted symmetric tree with daughter degree sequence (2,3,4) as shown in figure 4 below.

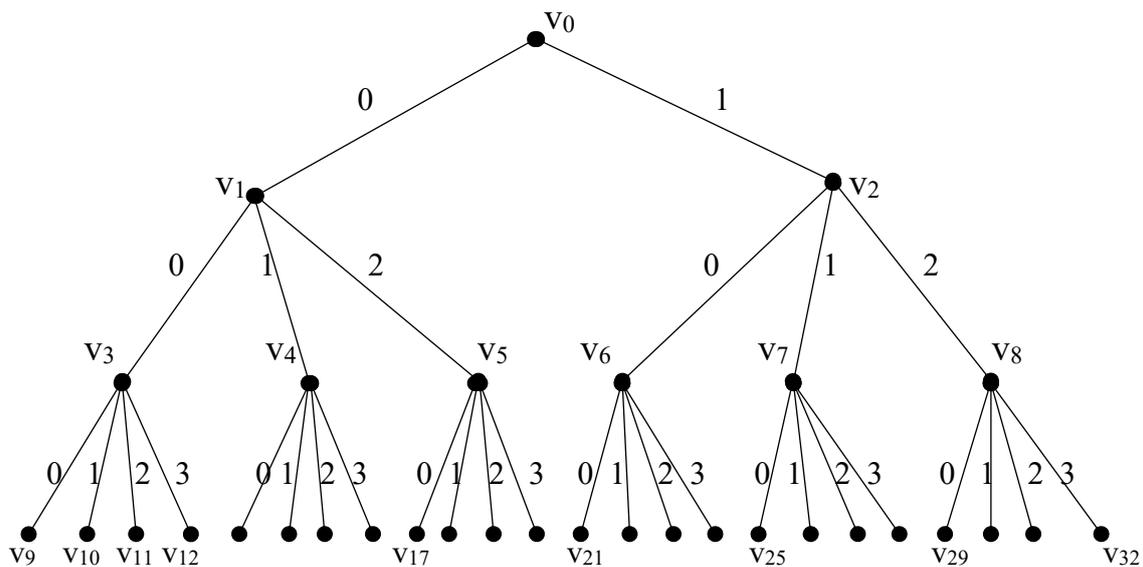

*Figure 4. (2,3,4) - rooted symmetric tree with indexed edges. Vertex $v_{21}$ in the bottom centre, for example, is identified by the edge number sequence (1,0,0)*

The following theorem provides a function that gracefully labels a rooted symmetric tree. Proof of the theorem can be found in the paper by the author (Rofa, 2021). This is followed by an example showing how the theorem is applied.

**Theorem 1.** Let T be a rooted symmetric tree as defined and edge indexed above and let v be a vertex at level r where $1 \leq r \leq q$. Then, the function *the f: V → {0, 1, ..., |E|} defined by*

1) $f(v) = (k_1 - x_1)h_2 - x_2h_3 - \ldots - x_{r-1}h_r - \frac{r-2}{2}$, *if r is even, and*

2) $f(v) = x_1h_2 + x_2h_3 + \ldots + x_{r-1}h_r + \frac{r-1}{2}$, *if r is odd.*

defines a graceful labelling of the rooted symmetric tree T. *Where, for $1 \leq i \leq q-1$, $h_i = 1 + k_i + k_ik_{i+1} + \ldots + k_ik_{i+1}\ldots k_{q-1}$ and $h_q = 1$. The numbers $h_1, h_2, \ldots, h_q$ are called the level numbers of T.* (Note that for any given vertex v, say at level r between 1 and q inclusive, $h_r$ is the number of vertices of the rooted symmetric subtree whose root vertex is v, (thus $h_1=|VT|$ and for r between 2 and q inclusive, $h_r = (h_{r-1} -1)/k_{r-1}$.

(Notice further that $f(v_0) = 0$ as root vertex $v_0$ is on level 1, that is r = 1 and hence equation 2 applies).

**Example 1.** Consider the rooted symmetric tree with daughter degree sequence (2,3,4) shown and edge indexed as in figure 3, by definition we have

$h_1 = 1 + k_1 + k_1k_2 + k_1k_2k_3 = 1 + 2 + 2 \times 3 + 2 \times 3 \times 4 = 33$

$h_2 = 1 + k_2 + k_2k_3 = 1 + 3 + 3 \times 4 = 16$

$h_3 = 1 + k_3 = 1 + 4 = 5$

$h_4 = 1$

Applying the labelling scheme of Theorem 1 results in the graceful labelling that is shown in figure 5 below. For example, vertex $v_{32}$ which has the edge index array (1,2,3) at level 4 will be assigned the label given by equation 2 of Theorem 1 which is:

$(k_1 - x_1)h_2 - x_2h_3 - x_3h_4 - \frac{4-2}{2} = (2 - 1).(16) - (2).(5) - (3).(1) - 1 = 2.$

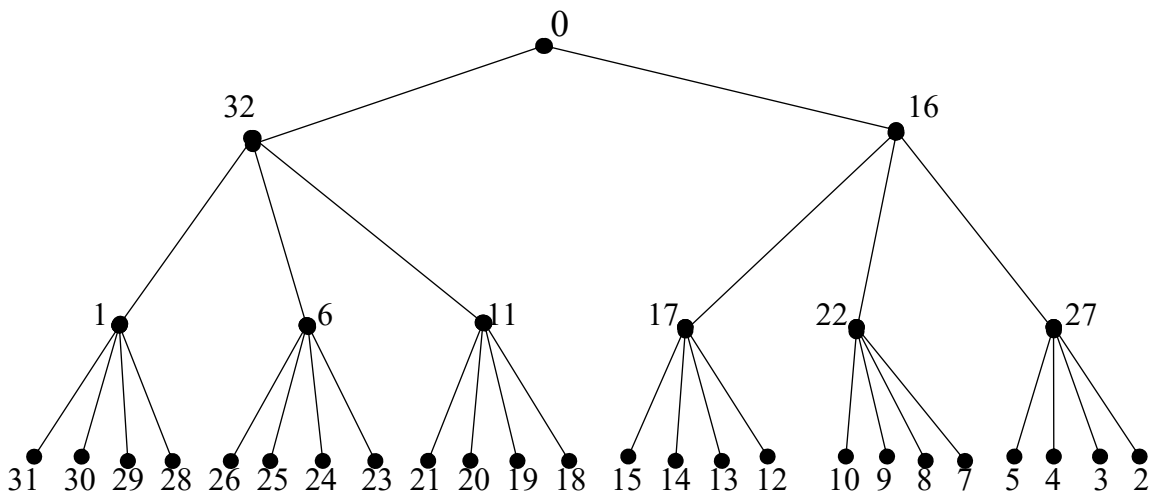

*Figure 5. A graceful labelling obtained by applying the function of Theorem 1 and based on the edge indices of figure 4*

Notably, a path with n vertices, $P_n$, is a rooted symmetric tree with daughter degree sequence (1,1,…,1) [repeated n-1 times] and it is interesting to note that application of the function f of Theorem 1 to $P_n$ produce the canonical graceful labelling of paths as described algorithmically by Rosa (Rosa, 1967)(See figure 6). Partial motivation of this paper for investigating the 0-rotatability of rooted symmetric trees lies behind the fact that paths, as a special class of rooted symmetric trees, are 0-rotatable which was proved by Rosa in 1977, about 10 years after his original paper of 1967 (Rosa, 1977)

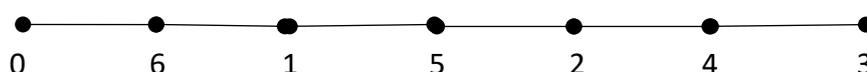

*Figure 6. Graceful labelling of $P_7$ is obtained either by canonical labelling by A. Rosa (Rosa, 1967) or through application of the algebraic function of Theorem 1 (Rofa, 2021)*

Figure 7 below shows various graceful labellings of $P_7$ with the label 0 assigned to various (and each) vertex of $P_7$.

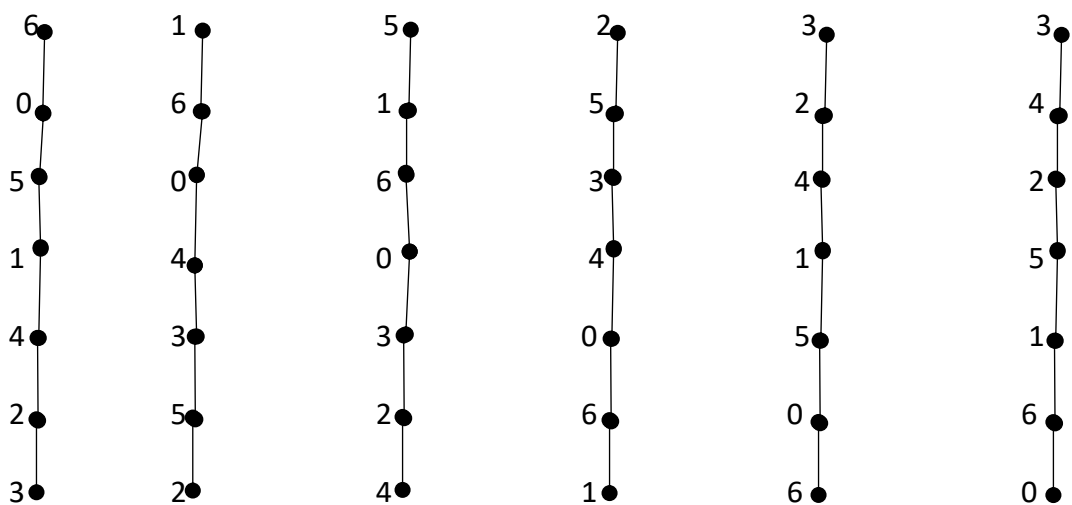

*Figure 7. Various graceful labellings of $P_7$ indicating its 0-rotatability*

## 3. Main Results

In this section we prove the following 3 objectives (as stated in the introduction):

(1) Rooted symmetric trees with 3 levels, including root vertex (that is, diameter 4) are 0-rotatable. We also prove the following about subclasses of rooted symmetric trees of diameter 6 as our 2nd and 3d objectives:

(2) Symmetric spiders with leg length at most 3 are 0-rotatable.

(3) Symmetric banana trees are 0-rotatable.

We first prove the following general theorem that will be used in all subsequent results:

**Theorem 2.** Given a rooted symmetric tree T with daughter degree sequence $(k_1, k_2, …, k_{q-1})$. If the tree P obtained by removing a $(k_1-1, k_2, …, k_{q-1})$ - rooted symmetric subtree (except root vertex) from T is a caterpillar, then T has graceful labellings that assign 0 to any given vertex at its last or second last levels (that is level q and level q-1) and it has graceful labellings that assign n-1 to any vertex at these levels.

Proof.

Case (1): T has an odd number of levels, that is q is odd. Consider any caterpillar P obtained by removing a $(k_1-1, k_2, …, k_{q-1})$ - rooted symmetric subtree (except root vertex) from T. Assuming that $|VP| = p$, then the vertices of P can be canonically gracefully labelled by using a function g with the consecutive vertex labels $0, 1, …, k_{q-1}, …, k_{q-1} + (q-3)/2$ that are joined to consecutive vertex labels $n-1, n-2, …, n - (q-1)/2$ with label 0 assigned to vertex at level q and label $k_{q-1} + (q-3)/2$ assigned to root vertex $v_0$. This induces the consecutive edge labels $n-1, n-2, n - (k_{q-1} + q - 2)$. Now, $H = T/P \cup \{v_0\}$ is a rooted symmetric tree with daughter degree sequence $(k_1-1, k_2, …, k_{q-1})$, Therefore by theorem 1, there is a graceful labelling h of H that assigns 0 to $v_0$, and producing the edge labels $1, 2, …, n - (k_{q-1} + q - 1)$. Now, increase every label of h by $k_{q-1} + (q - 3)/2$ to get a labelling h' of H which preserves the edge labels of h. It is easy to verify that increasing each vertex label of h by $k_{q-1}+(q-3)/2$ results in vertex labels that are disjoint from the vertex labels of g and complement them to get all the vertex labels from 0 to n-1 inclusive. [Minimum increased label (apart from $0 + k_{q-1}+(q-3)/2$) is $1+ k_{q-1}+(q-3)/2$ which is greater than the greatest vertex label of g and the maximal increased vertex label is less than $n-(k_{q-1}+q-2) + k_{q-1}+(q-3)/2 = n-(q-1)/2$] Hence extending the labelling h' of H with the labelling g of P gives a graceful labelling of T that assigns the label 0 to vertex v at level q and the label n-1 to an adjacent vertex at level q-1. A graceful labelling that interchanges the vertex labels 0 and n-1 is obtained by replacing each vertex label b of f by n-1-b. Now the result follows from the fact that all $(1, k_2, …, k_{q-1})$- rooted symmetric subtrees of T are isomorphic.

(note that the label $k_q+(q-3)/2$ is assigned to vertex $v_0$).

<u>Case (2):</u> T has an even number of levels (including root vertex). Again, consider caterpillar P obtained by removing a $(k_1-1, k_2, …, k_{q-1})$ - rooted symmetric subtree (except root vertex) from T. Assuming that $|VP|= p$, then the vertices of P can be canonically gracefully labelled by using a function g with the consecutive vertex labels 0, 1, …, $k_{q-1}$,…,$k_{q-1} + (q - 4)/2$ each of which adjacent to some of the consecutive edge labels n-1, n-2, …, n - q/2, with label 0 assigned to a vertex at level q and label n - q/2 assigned to root vertex $v_0$ and thus inducing the consecutive edge labels n-1, n-2, …, n - (q + $k_{q-1}$ - 2). Now, H = T/P ∪{$v_0$} is a rooted symmetric tree with daughter degree sequence $(k_1-1, k_2, …, k_{q-1})$, Therefore by theorem 1, there is a graceful labelling h of H that assigns 0 to $v_0$, and producing the edge labels 1, 2, …, n - (q + $k_{q-1}$ - 1). By replacing every label b of h by n – q/2 - b we get a labelling h' of H which preserves edge labels of h and which when extended with the labelling g of P gives a graceful labelling f of T that assigns the label 0 to vertex v at level q and n-1 to an adjacent vertex at level q-1. (Note that a similar argument as that presented in case (a) can be used to show that the vertex labels of h' complement those of g, except that root vertex has the same label in both, to form the set of all edge labels of T). A graceful labelling that interchanges the vertex labels 0 and n-1 is obtained by replacing each vertex label b of f by n-1-b. Now the result follows from the fact that all $(1, k_2, …, k_{q-1})$- rooted symmetric subtrees of T are isomorphic.
End of proof.

**Corollary 1.** Rooted symmetric trees with 3 levels (including root vertex as the first level) are 0-rotatable.

Proof. let T be a rooted symmetric tree on n vertices with root vertex $v_0$ and daughter degree sequence $(k_1, k_2)$. Theorem 1 provides a function f, that gracefully labels T with label 0 assigned to the root vertex $v_0$. Furthermore, the tree P obtained by removing a $(k_1-1, k_2)$- rooted symmetric subtree (except root vertex) from T is a caterpillar. Therefore, given a vertex at level 2 or level 3, then by Theorem 2 (case (1) of the proof), there exists a graceful labelling that assigns 0 to v. Therefore, T is 0-rotatable. This proves objective 1 as set out in the introduction. We illustrate Corollary 1 (and part (1) of Theorem 2) in the following example.

**Example 2.** Consider the rooted symmetric tree T with daughter degree sequence (3,4) [3 levels] that is shown in figure 8(a) below. By Theorem 1, there exists a graceful labelling that assigns 0 to root vertex. To obtain a graceful labelling that assigns 0 to a vertex at level 3, we do the following (as prescribed by case (1) of Theorem 2):

**Step 1**: Decompose the given (3,4) - rooted symmetric tree into a (2,4) - rooted symmetric tree (figure 8(b)) and a caterpillar (figure 8(c)).

**Step (2):** Gracefully label the caterpillar using canonical labelling of caterpillars as shown in figure 8(c) and gracefully label the (2,4) - rooted symmetric tree by using the labelling scheme of Theorem 1 as show figure 8(b).

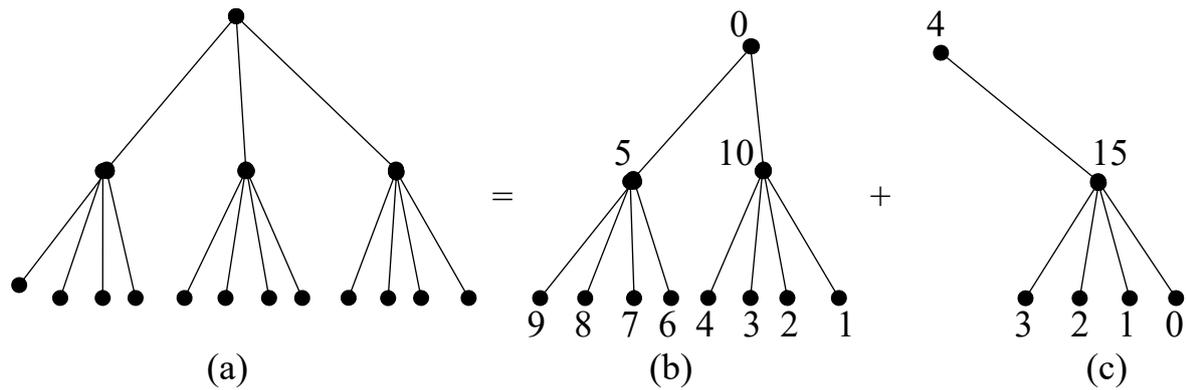

*Figure 8(a). (3,4) rooted symmetric tree decomposed into (2,3) – rooted symmetric tree(b) and caterpillar(c).*

**Step (3):** Increasing each vertex label of the tree in figure 8(b) by 4 and join the trees of (b) and (c) of figure 8 by identifying the root vertex to get the graceful labelling of the given (3,4) rooted symmetric tree which assigns the label 0 to a vertex at level 3, as shown in figure 9 below. A graceful labelling that assigns 0 to a vertex at level 2 is obtained by replacing each vertex label *b* by *15-b*. Since all (1,3) – rooted symmetric subtrees of T are isomorphic, it follows that T is 0-rotatable.

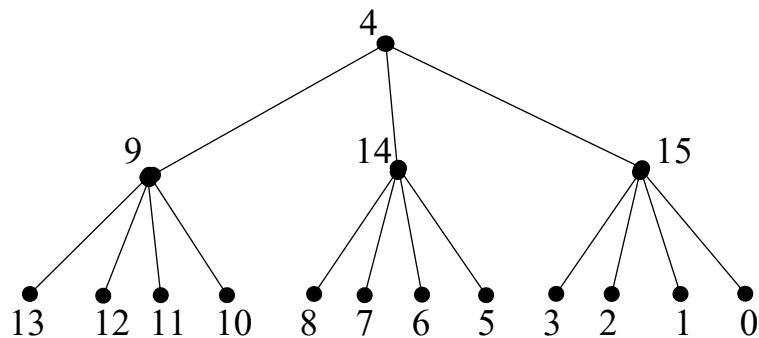

*Figure 9. A graceful labelling of the (4,3) rooted symmetric tree that assigns 0 to a vertex at level 3.*

We present here, for the record, an alternative method that asserts corollary 1 but without proof. The reason it is presented here is that in addition to asserting Corollary 1, it could be useful towards proving or disproving our conjecture that rooted symmetric trees are 0-rotatable.

**Lemma 1**. let T be a rooted symmetric tree with 3 levels with daughter degree sequence ($k_1$, $k_2$) and level numbers $h_1$, $h_2$, and $h_3$. Let f be the graceful labelling of T according to the labelling prescribed by Theorem 1 that assigns 0 root vertex, then the product of the transpositions

$$(0 \quad k_2)(h_2 \quad h_2+k_2)(2h_2 \quad 2h_2+k_2) \ldots ([k_1-1]h_2 \quad [k_1-1]h_2+k_2)$$

is a permutation of the vertex labels of f that yields a graceful labelling of T that assigns the label 0 to a vertex at level 3.

We illustrate Lemma 1 by the next example which uses the same tree of example 2.
**Example 3.** For the (3,4)-rooted symmetric tree that was used for example 2 above: $k_1=3$, $k_2=4$,
$h_1= 1+3+12=15$,
$h_2=1+4=5$, and
$h_3=1$.
Application of Theorem 1 gives the graceful labelling *f* shown in 10(a) which assigns 0 to root vertex. Application of the permutation $(0 \ 4)(5 \ 9)(10 \ 14)$ of Lemma 1 to the vertex labels of *f* yields a graceful labelling which assigns the label 0 to a vertex at level 3 as shown in figure 10(b).

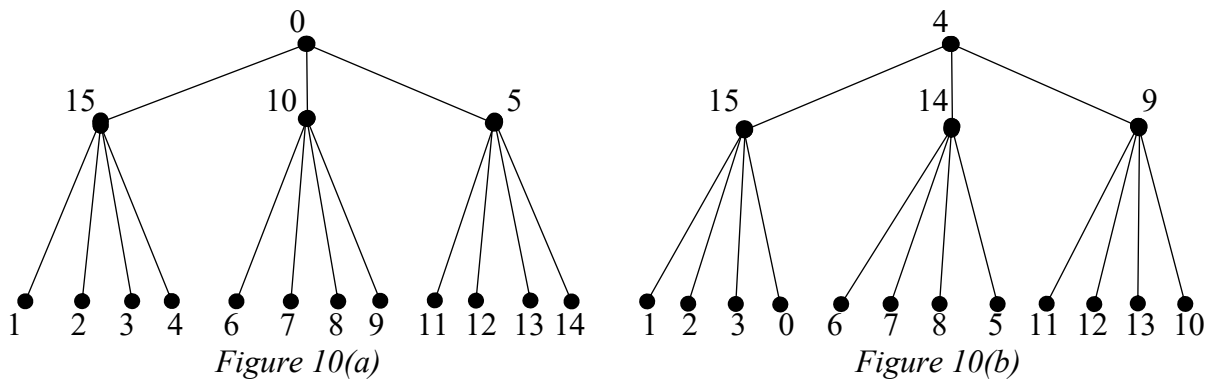

*Figure 10(a)*        *Figure 10(b)*

For the purpose of the next corollary, we have the following definition. Given a subset W of VT we say that *T is 0-rotatable within W* if for any v in W there exists a graceful labelling that assigns 0 to v. We have the following corollary.

**Corollary 2**. A rooted symmetric tree with daughter degree sequence $(k_1, 1, \ldots, 1, k_{q-1})$ is 0-rotatable within its first two and last two of its levels.

**Proof.** let T be a rooted symmetric tree on n vertices with root vertex $v_0$ and daughter degree sequence $(k_1, 1, \ldots, 1, k_{q-1})$. Theorem 1 provides a function $f$, that gracefully labels T, with label 0 assigned to the root vertex $v_0$. A graceful labelling that assigns 0 to an adjacent vertex at level 2 is obtained by replacing each vertex label $b$ of $f$ by $n - 1 - b$. Since all $(1, 1, \ldots, 1, k_{q-1})$ - rooted symmetric subtrees of T are isomorphic it follows that T is 0-rotatable within its first two levels. Furthermore, the tree P obtained by removing a $(k_1-1, 1, \ldots, 1, k_{q-1})$ rooted symmetric subtree (except root vertex) from T is a caterpillar therefore, by Theorem 2, T is 0-rotatable within its last two levels which completes the proof.

Next, we state and prove objective 2 as stated in the introduction.

**Corollary 3.** Symmetric Spiders with leg length at most 3 are 0-rotatable

**Proof.** A symmetric spider with leg length 1 or 2 is a rooted symmetric tree with at most 3 levels which, by Corollary1, is 0-rotatable. A symmetric spider with leg length 3 is a rooted symmetric tree with 4 levels satisfying the conditions of Corollary 2 and hence is 0 - rotatable. We illustrate Corollary 3 in the following example.

**Example 4.** Figure 11 below shows a symmetric spider with leg length 3 and with 4 labellings which are shown in figure 11(a): label 0 is assigned to root vertex (By Theorem 1), figure 11(b): label 0 is assigned to a vertex at level 2, figure 11(c): label 0 is assigned to a vertex at level 4 and figure 11(d): label 0 is assigned to a vertex at level 3.

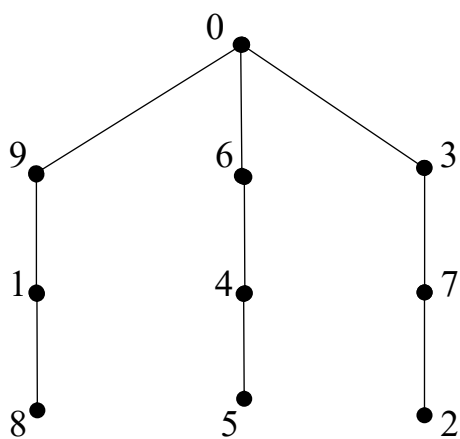

Figure 11(a)

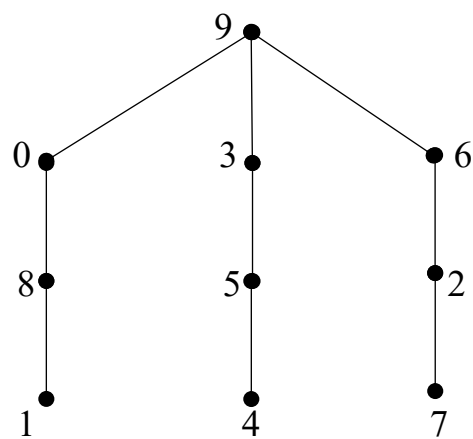

Figure 11(b)

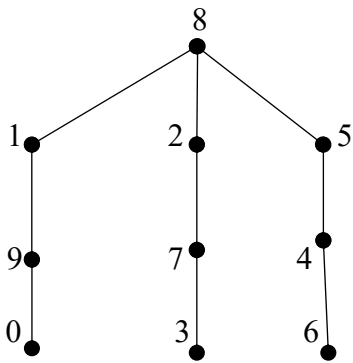
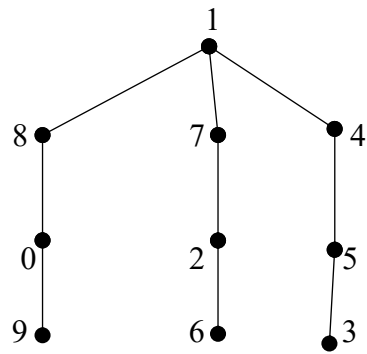

Figure 11(c)                    Figure 11(d)

Next, we state and prove objective 3 as stated in the introduction.

**Corollary 4.** Symmetric banana trees are 0-rotatable.

Proof. A symmetric banana tree is a $(k_1, k_2 = 1, k_3)$ - rooted symmetric tree for some positive integers $k_1$ and $k_3$, Hence, it is 0-rotatable by Corollary 2.

We illustrate Corollary 4 in the following example.

**Example 5.** Figure 12 below shows various graceful labellings of a banana tree with label 0 assigned to a vertex at each level.

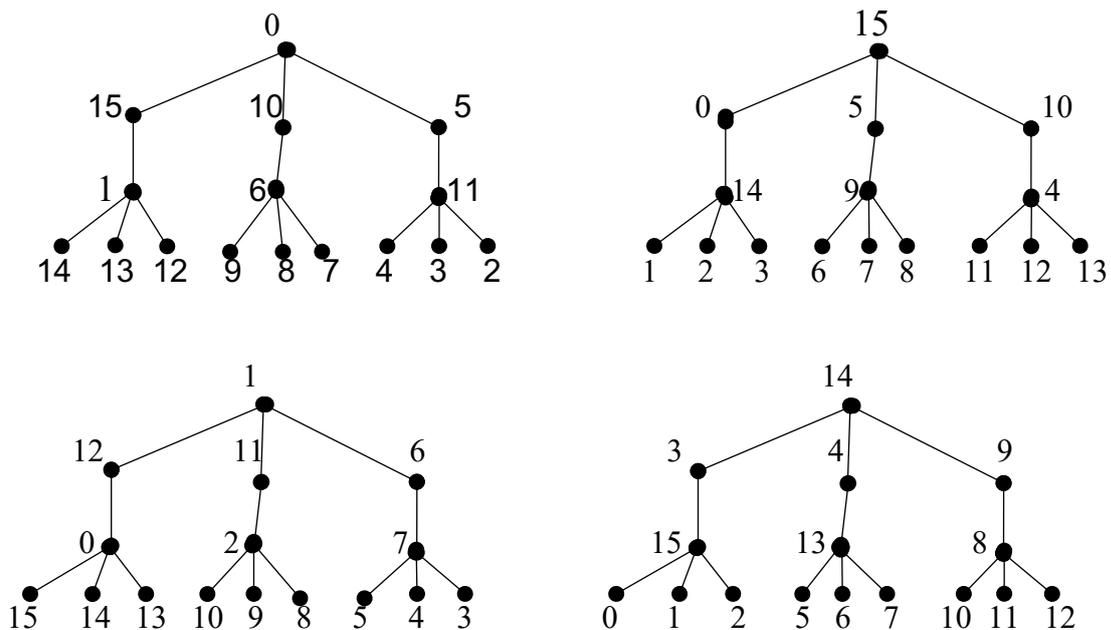

*Figure 12. Graceful labellings of a banana tree with label 0 assigned to a vertex at each level*

## 4. Summary

In this paper we have proved the following
- rooted symmetric trees with 3 levels are 0-rotatable.

We also proved the following about subclasses of rooted symmetric trees
- Symmetric spider trees with leg length at most 3.
- Symmetric banana trees.

Based on these results and within the context that rooted symmetric trees are generalisations of paths which have been shown by A. Rosa to be 0-rotatable, we propose the overarching conjecture that all rooted symmetric trees are 0-rotatable. We also propose the milestone sub-conjecture that all symmetric spider trees are 0-rotatable and since it has been shown in this paper that symmetric spider trees of leg length at most 3 are 0-rotatable, it remains to investigate symmetric spider trees of leg length 4 or more in order to settle the mentioned sub-conjecture.